\documentclass[journal]{IEEEtran}
\IEEEoverridecommandlockouts

\usepackage{calc}
\usepackage{bm}
\usepackage{amsmath}
\usepackage{amssymb}
\usepackage{amsthm}
\usepackage{amsfonts}
\usepackage{multirow}		
\usepackage{stfloats}
\usepackage{xfrac}			
\usepackage{booktabs}

\usepackage{algpseudocode}
\usepackage{algorithm}

\usepackage{color,soul}

\newcommand{\na}{\nabla}

\usepackage{mathtools}
\usepackage{cuted}
\usepackage{url}	
\usepackage{accents}

\allowdisplaybreaks

\usepackage{wrapfig}

\usepackage[noadjust]{cite}

\PassOptionsToPackage{hyphens}{url}
\usepackage{hyperref}

\begin{document}
\vspace{-30pt}

\title{Smoothed Two-Stage Decomposition Algorithm for Solving Large-Scale Transmission and Distribution AC-OPF Problems}

\author{\vspace{-5pt}Juan Ospina, Manuel Garcia, Xinyi Luo, Andreas W\"{a}chter, David M. Fobes, Russell Bent \vspace{-30pt}
\thanks{
(\textit{Corresponding author}: Juan Ospina )

J. Ospina, M. Garcia, D. Fobes, and R. Bent are with Los Alamos National Laboratory, Los Alamos, NM 87545-1663, USA (email: jjospina@lanl.gov; mjgarcia@lanl.gov; dfobes@lanl.gov; rbent@lanl.gov)

X. Luo and A. W\"achter are with the Northwestern University, Evanston, IL 60208, USA (email: xinyiluo2023@u.northwestern.edu; andreas.waechter@northwestern.edu)

}

}

\maketitle

\begin{abstract}
The integration of distributed energy resources (DERs) into the power grid has introduced new challenges to AC optimal power flow (AC-OPF) problems. Traditional OPF optimize consider transmission systems, treating distribution networks as static loads. However, the growing presence of DERs makes accurate distribution system modeling crucial for grid operations. Consequently, efficiently solving the resultant large-scale, nonconvex transmission and distribution (T\&D) AC-OPF problem remains a significant challenge. This paper proposes a Smoothed Two-Stage Decomposition Optimizer (StsDOpt) to address these complexities by decomposing the T\&D AC-OPF problem into a master-subproblem(s) structure, enabling parallel solving. Unlike traditional methods, StsDOpt does not rely on approximations or relaxations. It uses a smoothing technique to render the subproblems’ responses differentiable with respect to the master problem, leveraging the barrier problem properties inherent in primal-dual interior point methods. This approach is crucial for accurately modeling and solving distribution systems, which are multiphase, unbalanced, and nonlinear, distinguishing StsDOpt apart from other methods. Integrated into the PowerModelsITD framework, StsDOpt has been validated through numerical experiments, demonstrating reduced wall-clock solve time and increased scalability. Results highlight its efficacy as a robust, scalable solution for large-scale T\&D AC-OPF problems, facilitating the reliable integration of DERs into complex T\&D systems.
\end{abstract}


\vspace{-6mm}
\section{Introduction}\label{Sec:1}

The integration of distributed energy resources (DERs) into the power grid is essential for meeting rising energy demands and advancing decarbonization, energy affordability, energy efficiency, and resiliency goals. Energy demands have been rising in the last couple of years, driven by population growth, increased data center demands, and electric vehicle (EV) adoption, among others. Estimates made by various agencies in 2023 indicate that electric load in the United States is expected to grow 4.7\% over the next five years; an estimate that was corrected from the 2.6\% made in 2022 \cite{utilitydive, eraflatpower}. Longer-term estimates forecast that, in the US, load is expected to grow 42\% by 2050 \cite{enverus}. Many experts advise that the electric grid is `not prepared' for this significant load growth and that `low transfer capabilities between regions' are a key risk for reliability if new generation resources are not integrated into the grid \cite{utilitydive}. That is why, initiatives, such as the Inflation Reduction Act and the Federal Energy Regulatory Commission’s  (FERC) Order No. 2222 \cite{ferc2222}, propose policy and market structures that enable DERs to provide services to the grid in exchange for financial compensation with the objective of incentivizing their rapid adoption for alleviating transfer capability requirements between regions and creating a value stream for customers and entities that use and deploy these resources \cite{deloitte, wri}.

Nonetheless, as DERs become more prevalent, distribution systems are becoming more complex and difficult to operate since their operation is becoming intertwined with the operation of the transmission system. Traditionally, transmission systems are overseen by a system operator, commonly known as a Transmission System Operator (TSO), that consider the distribution system to be a passive energy requester composed of inelastic demand that is forecasted and treated as being fixed. However, as distribution systems become more `active' due to the introduction of DERs, the traditional assumption of passive distribution systems is no longer valid and transmission and distribution (T\&D) co-optimization methodologies become essential for coordinating and optimizing the use of generation resources across T\&D operational boundaries.

The AC optimal power flow (AC-OPF) problem, which underpins many critical transmission operational functions such as Unit Commitment and Real-time Economic Dispatch, presents significant computational challenges for large-scale systems due to its inherent nonconvexity and limited computationally tractability. In such large-scale AC-OPF problems, computational challenges are further exacerbated when multi-phase unbalance distribution systems are included to create large-scale T\&D AC-OPF problems that enable the coordination and optimization of resources at both the transmission and distribution system levels. As a result, traditional solvers and methods struggle to find efficient or even feasible solutions within reasonable computational times, hindering their ability to perform cutting-edge T\&D operational functions.   
\vspace{-2mm}
\subsection{Literature Review}

The computational challenges encountered when solving large-scale T\&D AC-OPF problems have motivated many researchers to design novel computational methods and algorithms capable of handling large-scale T\&D models. The most common approach explored by researchers decomposes the T\&D AC-OPF problem into a transmission-level master OPF problem and distribution-level subproblems that interact with each other. For instance, authors in \cite{decentralizedacoptimal} propose a decentralized distribution-cost correction methodology that solves the decomposed AC-OPF using an AC formulation for the transmission system master problem and second-order cone (SOC) branch flow (BF) approximations to represent the distribution system subproblems. Authors in \cite{coordinatedTanD} use a heterogeneous decomposition algorithm while researchers in \cite{sun2014master} propose an approach that solves the OPF problem using a decomposition and splitting-based iterative method, where admittance equivalencing is used for modeling the distribution systems connected to the transmission system.  

Other approaches, such as the one presented in \cite{yuan2017hierarchical}, utilize Benders decomposition to solve a convex approximation of the AC-OPF problem based on the SOC formulation. A similar approach is presented in \cite{NAWAZ2024110005}, where researchers use Benders decomposition to solve a risk-based stochastic security-constrained unit commitment problem, dividing master and subproblems to generate feasibility cuts. A decomposition approach based on Kron reductions and Benders-type decomposition techniques is presented in \cite{CONSTANTEFLORES2024427}. This approach is also used to solve a security-constrained unit commitment problem.

The alternative approach is to solve the T\&D AC-OPF problem as an integrated problem where all transmission and distribution system models are combined into a monolithic optimization model solved. The primary example of this approach is presented in \cite{ospina2024modeling}, where an integrated T\&D power network optimization framework is proposed. This optimization framework enables the computational evaluation of a diverse set of power network formulations (e.g., AC polar, AC rectangular, Current-Voltage) and problem specifications (e.g., OPF, PF) using a common platform. This research demonstrated that solving T\&D AC-OPF problems using monolithic models is not a significant challenge until reaching 50k+ nodes or 1M+ variables. The PowerModelsITD \cite{ospina2024modeling} optimization framework serves as the baseline used to evaluate the performance and scalability of the smoothed two-stage decomposition optimizer (StsDOpt) proposed in this research. Other researchers have also explored the use of integrated T\&D models for stochastic optimization \cite{10517420, 10574270} and planning \cite{10089624}. A comprehensive review of the challenges behind three T\&D coordination strategies, i.e., transmission-managed, distribution-managed, and transmission-distribution hybrid, is presented in \cite{givisiez2020review}.

\vspace{-2mm}
\subsection{Contributions}

While existing methods for solving T\&D OPF problems have made significant contributions to the field, our proposed decomposition-based optimizer (i.e., StsDOpt) offers several distinct advantages when compared to current literature for solving large-scale nonlinear nonconvex T\&D AC-OPF problems. Unlike most reviewed research that rely on approximations and/or relaxations, we propose a method that can handle nonlinear nonconvex models when decomposing the AC-OPF problem, i.e., there is no need for linearization or relaxation, thus there is no need for checking AC feasibility; a required time-consuming post-processing step needed when using approximations and/or relaxations. This is also particularly important for accurately modeling distribution systems, which are inherently unbalanced and nonlinear. Furthermore, while many of the research papers reviewed assume single-phase balanced models when modeling distribution systems, our proposed method can accommodate realistic multi-phase unbalanced systems. These contrasting characteristics also differentiate the proposed StsDOpt from the foundational work done in \cite{tu2020two}. Note that a direct numerical comparison with every existing method is not feasible and beyond the scope of this work; as they differ in their fundamental representation of the T\&D models, e.g., using single-phase representations. Instead, this paper focuses on establishing theoretical properties and demonstrating scalability in a series of large‑scale test cases. Here is a list of specific contributions of this research and differences from previous work.

\vspace{-2mm}
\begin{enumerate}
    \item The proposed approach directly solves nonlinear, nonconvex optimization models without approximations or relaxations, unlike existing methods that typically assume single-phase balanced systems. This enables accurate modeling of realistic multi-phase unbalanced distribution networks.
    
    \item By avoiding approximations and relaxations, the proposed method eliminates the need for computationally expensive AC feasibility checks.
    
    \item A Sequential Quadratic Programming (SQP) solver \cite{xinyi_thesis} is developed for the smoothed master problem and tightly integrated with the interior-point optimizer Ipopt \cite{wachter2006implementation} for solving the subproblems. Parallelization is implemented using Julia. The proposed approach significantly reduces wall-clock time compared to monolithic solution approaches.
    
    \item The implementation exploits Ipopt's efficient internal linear algebra routines to compute derivatives, reducing computational complexity while avoiding nonconvexity issues reported in previous work that can slow convergence of the master problem.
    
    \item Unlike \cite{tu2020two}, the proposed method requires only second-order derivatives rather than third-order derivatives, substantially reducing computational overhead.
    
    \item Extrapolation techniques are incorporated to provide improved initial points for the master problem, accelerating local convergence.
\end{enumerate}

In addition to these main contributions, our proposed method also offers several additional key advantages over existing solutions. Firstly, the proposed StsDOpt solver is integrated into the PowerModelsITD optimization framework \cite{ospina2024modeling}, allowing for thorough testing and verification of the results. This is in contrast to many reviewed papers that do not provide effective ways to verify and replicate the presented results. Additionally, our approach is scalable and is capable of handling large-scale T\&D models, even when memory is scarce, by providing the ability to \textit{multithread} or \textit{distribute} the solver to improve performance. Moreover, the proposed decomposition-based method, as part of the PowerModelITD optimization framework, is modular and allows researchers to develop and deploy new network formulations (nonlinear, approximations, and/or relaxations) that can be seamlessly integrated into the optimization framework and can be solved using the StsDOpt optimizer.   

The remainder of the paper is organized as follows. Section \ref{sect:sts} presents the mathematical formulation of the proposed smoothed two-stage decomposition algorithm, the basis for the StsDOpt solver. Section \ref{sect:implementation} describes the implementation of the StsDOpt solver so that it can be interfaced and used within the Julia JuMP and PowerModelsITD environments. Section \ref{sect:Results} describes the experimental setup and the numerical and scalability results of the test cases used to evaluate the performance of the proposed method. Finally, Section \ref{sect:conclusions} concludes the paper and discusses the future directions of the work. 
\vspace{-2mm}
\section{Smoothed Two-Stage Decomposition Algorithm}
\label{sect:sts}

\subsection{Smoothed Two-Stage Decomposition Algorithm}

The proposed decomposition algorithm is capable of solving nonlinear nonconvex continuous two-stage optimization problems where the first-stage (or master) problem has the form:

\begin{equation}\label{eq:first-stage}
    \begin{aligned}
        \min_x\quad & f_0(x) + \sum_{i=1}^N \hat f_i(x) \\
        {s.t.} \quad & g_0(x) = 0,\; h_0(x) \leq 0
    \end{aligned}
\end{equation}
with $N$ second-stage problems (or subproblems):
$i=1,\ldots,N$
\begin{equation}\label{eq:second-stage}
    \begin{aligned}
        \hat f_i(x) = \min_{y_1}\quad & f_i(y_i) \\
        {s.t.} \quad & g_i(x,y_i) = 0,\; h_i(x,y_i) \leq 0,
    \end{aligned}
\end{equation}
where $x$ and $y_i$ represent the variables of the optimization problem (master and subproblem, respectively), $f_i$ are the cost functions, $g_i$ the equality constraints, and $h_i$ the inequality constraints. 
It is important to note that the value functions $\hat f_i(x)$ are usually nonsmooth and not differentiable whenever the set of inequality constraints, that are active at the optimal solution of Eq.~\eqref{eq:second-stage}, changes with respect to $x$. This makes the use of efficient, second-order optimization methods for Eq.~\eqref{eq:first-stage} impossible. To overcome this issue, we follow a smoothing approach, introduced in \cite{tu2020two}, that replaces the second-stage problems with their \textit{barrier problem formulations} of the form
\begin{equation}\label{eq:second-stage-smoothed}
    \begin{aligned}
        \hat f_i(x;\mu) = \min_{y_i,s_i,\tilde x_i}\quad & f_i(y_i ) - \mu\sum_{j=1}^{n_i}\log(s_{ij}) \\
        {s.t.} \quad & g_i(\tilde x_i, y_i) = 0,\; h_i(\tilde x_i,y_i) + s_i = 0\\
        & x - \tilde x_i  = 0
    \end{aligned}
\end{equation}
with a fixed barrier parameter $\mu>0$ and slack variables $s_i$. Because there are no inequality constraints in this formulation, the optimal solutions of these problems depend \textit{smoothly} on the input variable vector $x$. Note that we introduced copies, $\tilde x_i$, of the first-stage variables because the multipliers related to the last constraint shown in Eq.~\eqref{eq:second-stage-smoothed} will help us compute derivatives of $\hat f_i$.

The barrier problem formulation, i.e.,  Eq.~\eqref{eq:second-stage-smoothed}, is at the heart of modern efficient interior-point methods for which powerful software implementations are readily available; e.g., Ipopt \cite{wachter2006implementation}.  It is well known, that, for a fixed input $x$, local minima of Eq.~\eqref{eq:second-stage-smoothed} converge to local minima of Eq.~\eqref{eq:second-stage} as the barrier parameters $\mu_k$ approaches zero, under standard assumptions. As a consequence, we can compute solutions of Eq.~\eqref{eq:first-stage} by solving a sequence of smoothed first-stage problems of the form:
\begin{equation}\label{eq:first-stage-smoothed}
    \begin{aligned}
        \min\quad & f_0(x) + \sum_{i=1}^N \hat f_i(x;\mu_k) \\
        {s.t.} \quad & g_0(x) = 0, h_0(x) \leq 0
    \end{aligned}
\end{equation}
for a decreasing sequence of barrier parameters $\mu_k\to0$. This approach has the computational advantage that efficient second-order algorithms can be utilized for both the first-stage, Eq.~\eqref{eq:first-stage-smoothed}, and second-stage, Eq.~\eqref{eq:second-stage-smoothed}, of the smoothed two-stage optimization problem.
In order to be able to apply these methods, first and second derivatives of $\hat f_i(x;\mu)$ must be provided. In the proposed decomposition approach, these are computed in the following manner.

When an interior-point solver is applied to a smoothed second-stage problem, Eq.~\eqref{eq:second-stage-smoothed} computes a primal-dual solution, $v_i=(y_i,s_i,\tilde x_i,\lambda_i^g ,\lambda_i^h,\eta_i)$, of the first-order optimality conditions of Eq. \eqref{eq:second-stage-smoothed},
\begin{multline*}
    F_i(x,v_i;\mu) = \\
    \begin{pmatrix}
        \na_{y_i} f_i(y_i) + \na_{y_i} g_i(x,y_i)\lambda_i^g  + \na_{y_i} h_i(x,y_i)\lambda_i^h\\
         \na_{x} g_i(x,y_i)\lambda_i^g  + \na_{x} h_i(x,y_i)\lambda_i^h + \eta_i\\
        g_i(\tilde x_i, y_i)\\
        h_i(\tilde x_i,y_i) + s_i \\
        x - \tilde x \\
        s_i\circ \lambda_i^h - \mu e
    \end{pmatrix}=0,
\end{multline*}
where $\lambda_i^g ,\lambda_i^h,\eta_i$ are the dual variables for the constraints in Eq. \eqref{eq:second-stage-smoothed}, $\circ$ denotes the component-wise product of two vectors, and $e$ represents the vector of ones. Standard sensitivity results imply that $\eta_i=\na_{x} \hat f_i(x;\mu)$.

For second derivatives, we can apply the implicit function theorem to the relationship $F_i(x,v_i(x);\mu)=0$ and compute $\na_x v_i(x)$ by solving the linear system:
\begin{equation}\label{eq:implicit-function}
\na_{v_i} F_i(x,v_i(x);\mu)^T\na_x v_i(x)^T = - \na_x F_i(x,v_i(x),\mu)^T.
\end{equation}
Then
$\na^2\hat f_i(x;\mu)$ can be extracted as the component $\na_x \eta_i(x)$ in $\na_x v_i(x)$.
A high-level description of the algorithm is given in Algorithm \ref{alg:decomposition_algo}.

\begin{algorithm}[t]
\caption{Smoothed two-stage decomposition algorithm} \label{alg:decomposition_algo}
\begin{algorithmic}[1]
\Require{Initial master problem iterate $x_0$; initial smoothing parameter $\mu_0$.}

\State Set $\tilde x_0\gets x_0$.
\For{$k=0,1,2,\ldots$}
\State Starting from $\tilde x_k$, (approximately) solve the smoothed master problem, Eq.~\eqref{eq:first-stage-smoothed}.
\label{line:master_solver}
\State Reduce $\mu_k$ to get $\mu_{k+1}<\mu_k$ (so that $\mu_k\to0$).\label{line:reduce_mu}
\State Choose starting iterate $\tilde x_{l+1}$ for the next iteration. 
\label{line:initialization_decay_mu}
\EndFor
\end{algorithmic}
\end{algorithm}

Algorithm \ref{alg:decomposition_algo} alternates between solving a smoothed subproblem (Stage 1) and solving the subproblems (Stage 2). Subproblems are solved independently for each transmission–distribution boundary, and their solutions feed directly into the master update.

In every inner iteration of the master problem solver in Step~\ref{line:master_solver}, the subproblems must be solved for the evaluation of $\na \hat f_i$ and its derivatives.
For every outer SQP iteration, it is possible to warm‑start the variables of the subproblem solver, significantly reducing the number of interior‑point iterations needed to reach the target tolerance.
Another important remark to consider is that the smoothed master problem needs to be solved to high accuracy only as $\mu_k$ converges, so  for large values of $\mu_k$, approximated solutions are sufficient \cite{xinyi_thesis}.

\vspace{-2mm}
\subsection{T\&D AC-OPF Formulation based on Smoothed Two-Stage Decomposition Algorithm}

In this subsection, the formulation of the T\&D AC-OPF based on the AC polar formulation is presented. Using the same structure shown in Eqs. \eqref{eq:first-stage} -- \eqref{eq:second-stage-smoothed}, the T\&D problem can be defined as: 

\begin{equation}\label{eq:masterACP}
    \begin{split}
        \min \quad & f_0(x^{\mathcal{T}}) + \sum_{i=1}^N \hat f_i(x^{\mathcal{T}}_{i,\beta^{^\mathcal{D}}})\\
        {s.t.} \quad &g_0(x^{\mathcal{T}}) = 0,
        \\
        &h_0(x^{\mathcal{T}}) \leq 0
    \end{split}
\end{equation}
with $N$ distribution system subproblems:
\begin{equation}\label{eq:subproblemACP}
    \begin{aligned}
        \hat f_i(x^{\mathcal{T}}_{i,\beta^{^\mathcal{D}}}) = \min_{y^{\mathcal{D}}_i}\quad & f_i(y^{\mathcal{D}}_i) \\
        {s.t.} \quad & g_i(y^{\mathcal{D}}_i) = 0,\;
        \\ \quad & h_i(y^{\mathcal{D}}_i) \leq 0\\
        & x^{{\mathcal{T}}}_{i,\beta^{^\mathcal{D}}} - y^{\mathcal{D}}_{i,\beta^{^\mathcal{T}}}  = 0,
    \end{aligned}
\end{equation}
where the superscript $^{\mathcal{T}}$ relates to the transmission system, $^{\mathcal{D}}$ to the distribution systems, and subscript $\beta$ indicates variables that belong to buses that exist in the T\&D boundary bus set. Note that the boundary-related variables $x^{\mathcal{T}}_{i,\beta^{^\mathcal{D}}}$ are part of the variable set $x^{\mathcal{T}}$ and the boundary-related variables $y^{\mathcal{D}}_{i,\beta^{^\mathcal{T}}}$ are part of the variable set $y^{\mathcal{D}}_i$.

In this formulation, $f_0(x^{\mathcal{T}})$ represents a transmission system cost function of the form: 
\begin{equation}\label{eq:transcost}
    \begin{split}
        f_0(x^{\mathcal{T}}) &= \min \bigg(\sum_{k \in G^{\mathcal{T}}} c_{2k}(P_{g,k}^{\mathcal{T}})^2 + c_{1k}(P_{g,k}^{\mathcal{T}}) + c_{0k} \bigg),
    \end{split}
\end{equation}
where $P_{g,k}$ represents the active power dispatch of generator $k$ in the set of transmission system generators $G^\mathcal{T}$, and $c_{ik}$ represent cost coefficients related to the specific generation unit $k$. 

Similarly, $f_i(y^{\mathcal{D}}_i)$, represents a distribution system cost function of the form:
\begin{equation}\label{eq:distrocost}
    \begin{split}
        \hat f_i(y^{\mathcal{D}}_i) = \min \bigg(\sum_{m \in G_i^{\mathcal{D}}} c_{1m}(\sum_{\varphi \in \Phi} P_{g,m}^{\mathcal{D},\varphi}) + c_{0m} \bigg),
    \end{split}
\end{equation}
where $P_{g,m}$ represents the active power dispatch of generator $m$ that belongs to the set of generators $G_i^\mathcal{D}$ in the $i$-th distribution system and $c_{im}$ represent cost coefficients related to the specific generation unit $m$, i.e., DER.  The set of phases in the multiphase unbalance distribution system is denoted $\Phi$ and each phase is represented by $\varphi \in \Phi$.

At the boundary of the transmission and distribution system problem, shared variables are passed from the master problem to the subproblems in order to formulate the last equality constraint of Eq. \eqref{eq:subproblemACP}. In the AC polar formulation, the variables shared at the boundary are the active power, reactive power, voltage magnitude, and voltage angle at the boundary buses, $\beta^\mathcal{T}$ and $\beta^\mathcal{D}$, respectively. For the transmission system, the shared variables are:

\begin{equation} \label{eq:sharedtrans}
\bigg[ P_{i,\beta^{^\mathcal{D}}}^{\mathcal{T}}, Q_{i,\beta^{^\mathcal{D}}}^{\mathcal{T}}, V_{i,\beta^{^\mathcal{D}}}^{\mathcal{T}}, \theta_{i,\beta^{^\mathcal{D}}}^{\mathcal{T}} \bigg] 
\end{equation}

where all these boundary (shared) variables are represented by $x^{\mathcal{T}}_{i,\beta^{^\mathcal{D}}}$. Similarly, for distribution system $i$, the shared variables can be defined as:

\begin{equation} \label{eq:shareddist}
\bigg[ p_{i,\beta^{^\mathcal{T}}}^{\mathcal{D}}, q_{i,\beta^{^\mathcal{T}}}^{\mathcal{D}}, v_{i,\beta^{^\mathcal{T}}}^{\mathcal{D}}, \theta_{i,\beta^{^\mathcal{T}}}^{\mathcal{D}} \bigg] 
\end{equation}

where all these boundary (shared) variables are represented by $y^{\mathcal{D}}_{i,\beta^{^\mathcal{T}}}$. For each distribution system $i$, the shared variables are integrated into the optimization model via the following constraints.

\begin{gather}
    \label{eq:acpboundaryrealpowerequality}
    \sum_{\varphi \in \Phi} p_{i,\beta^{^\mathcal{D}} }^{\mathcal{D},\varphi} +  p_{i,\beta^{^\mathcal{T}}}^{\mathcal{D}} = 0, \ \forall (\beta^{^\mathcal{T}},\beta^{^\mathcal{D}}) \in {\Lambda}\\
    \label{eq:acpboundaryreactivepowerequality}
    \sum_{\varphi \in \Phi} q_{i,\beta^{^\mathcal{D}}}^{\mathcal{D},\varphi} +  q_{i,\beta^{^\mathcal{T}}}^{\mathcal{D}} = 0, \ \forall (\beta^{^\mathcal{T}},\beta^{^\mathcal{D}}) \in {\Lambda} \ \ \\
    \label{eq:acpboundaryvoltagemagnitudeAequality}
    v_{i,\beta^{^\mathcal{T}}}^{\mathcal{D}} = v_{i,\beta^{^\mathcal{D}}}^{^{\varphi}}, \ \forall (\beta^{^\mathcal{T}},\beta^{^\mathcal{D}}) \in {\Lambda}, \varphi\in\Phi \ \\
    \label{eq:acpboundaryvoltageangleAequality}
    \theta_{i,\beta^{^\mathcal{T}}}^{\mathcal{D}} = \theta_{i,\beta^{^\mathcal{D}}}^{^{a}}, \ \forall (\beta^{^\mathcal{T}},\beta^{^\mathcal{D}}) \in {\Lambda} \ \\
    \label{eq:acpboundaryvoltageangleBequality}
    \theta_{i,\beta^{^\mathcal{D}}}^{^{b}} = (\theta_{i,\beta^{^\mathcal{D}}}^{^{a}} -120^{\circ}),  \ \forall \beta^{^\mathcal{D}} \in B^{^\mathcal{B}} \cap  B^{^\mathcal{D}} \\
    \label{eq:acpboundaryvoltageangleCequality}
    \ \theta_{i,\beta^{^\mathcal{D}}}^{^{c}} = (\theta_{i,\beta^{^\mathcal{D}}}^{^{a}} +120^{\circ}), \ \forall \beta^{^\mathcal{D}} \in B^{^\mathcal{B}} \cap  B^{^\mathcal{D}}
\end{gather}

\noindent Eqs.~\eqref{eq:acpboundaryrealpowerequality}--\eqref{eq:acpboundaryreactivepowerequality} equalize the active and reactive power flow values at the boundary, where the power flow from the distribution system is calculated as the sum over all phases (i.e.,  phases $\Phi \in \{a,b,c\}$) and $\Lambda$ is the set of boundary links in the overall optimization problem. 
Eqs.~\eqref{eq:acpboundaryvoltagemagnitudeAequality}--\eqref{eq:acpboundaryvoltageangleCequality} define the equality constraints related to the voltage magnitudes and angles of the boundary buses. These constraints allow the one-to-one mapping of the single-phase-to-three-phase boundary buses. $B^{^\mathcal{B}}$ is the set of buses that belong to boundary buses set and $B^{^\mathcal{D}}$ is the set of buses or nodes that belong to the respective distribution system. To avoid repetition, all other equality and inequality constraints (i.e., $g(\,\cdot\,)$ and $h(\,\cdot\,)$) related to the specific AC polar OPF formulation are described in detail in \cite{ospina2024modeling}.

To ensure that the subproblems are always feasible (e.g., when the power provided by the transmission system, combined with the power generated by DERs in the distribution system, are not sufficient to cover demand), the active and reactive power of the transmission and distribution systems are allowed to be different, i.e., we relax the constraints in Eq.~\eqref{eq:acpboundaryrealpowerequality} and \eqref{eq:acpboundaryreactivepowerequality}, but their difference is penalized using a penalty term that is added to the subproblem cost function.  When the corresponding {penalty factors} are chosen to be sufficiently large, the power flows at the transmission/distribution system boundary match at the optimal solution.  However, if the factors chosen are too large, they may cause numerical difficulties in the algorithms.

\section{Implementation of the Smoothed Two-Stage Decomposition Algorithm}

\label{sect:implementation}

\begin{figure}[ht]
\centering
\includegraphics[width = 0.48\textwidth]{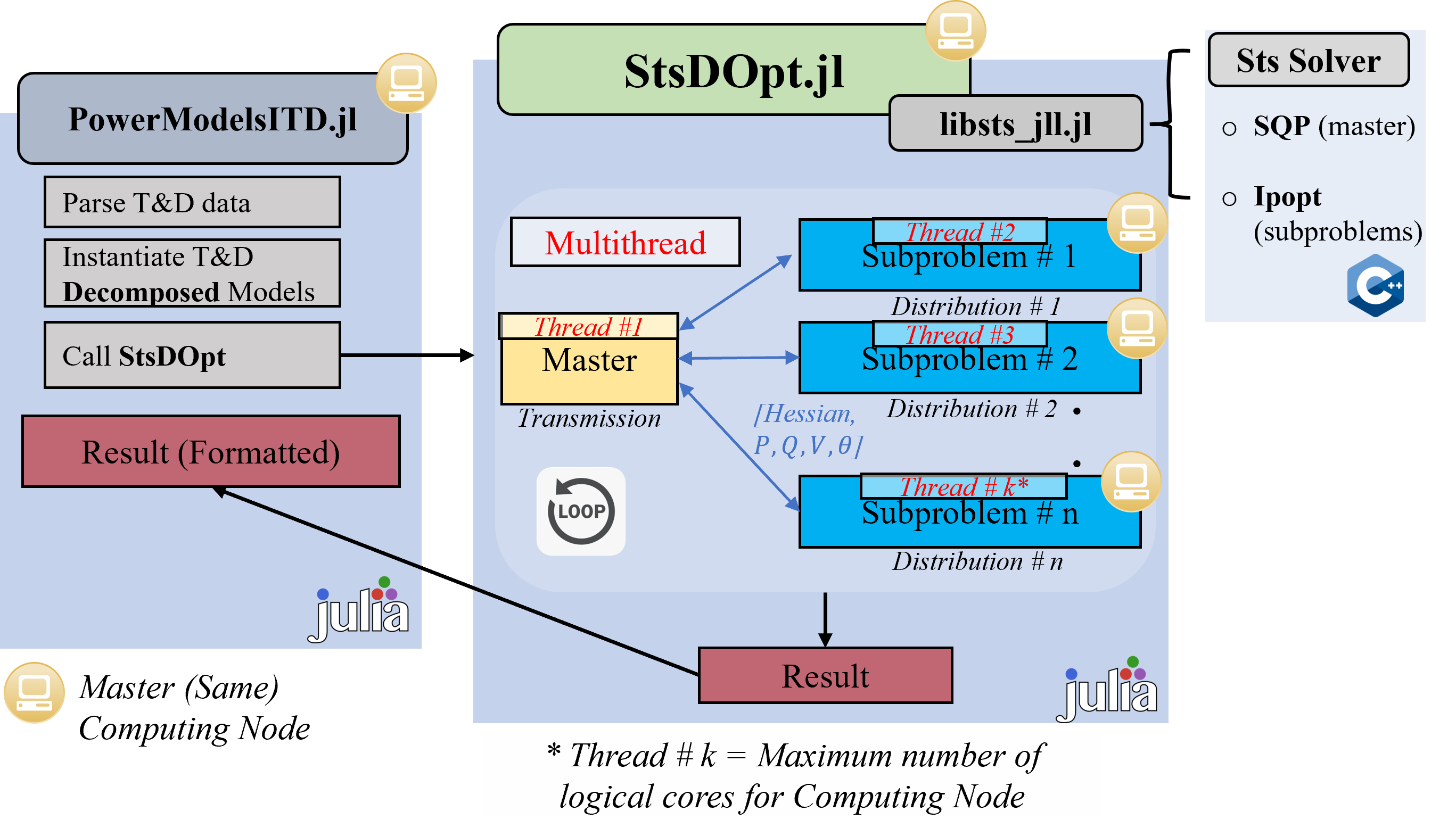}
\caption{\label{fig:multithread} Multithread version of StsDOpt.}
\end{figure}

\begin{figure}[ht]
\centering
\includegraphics[width = 0.48\textwidth]{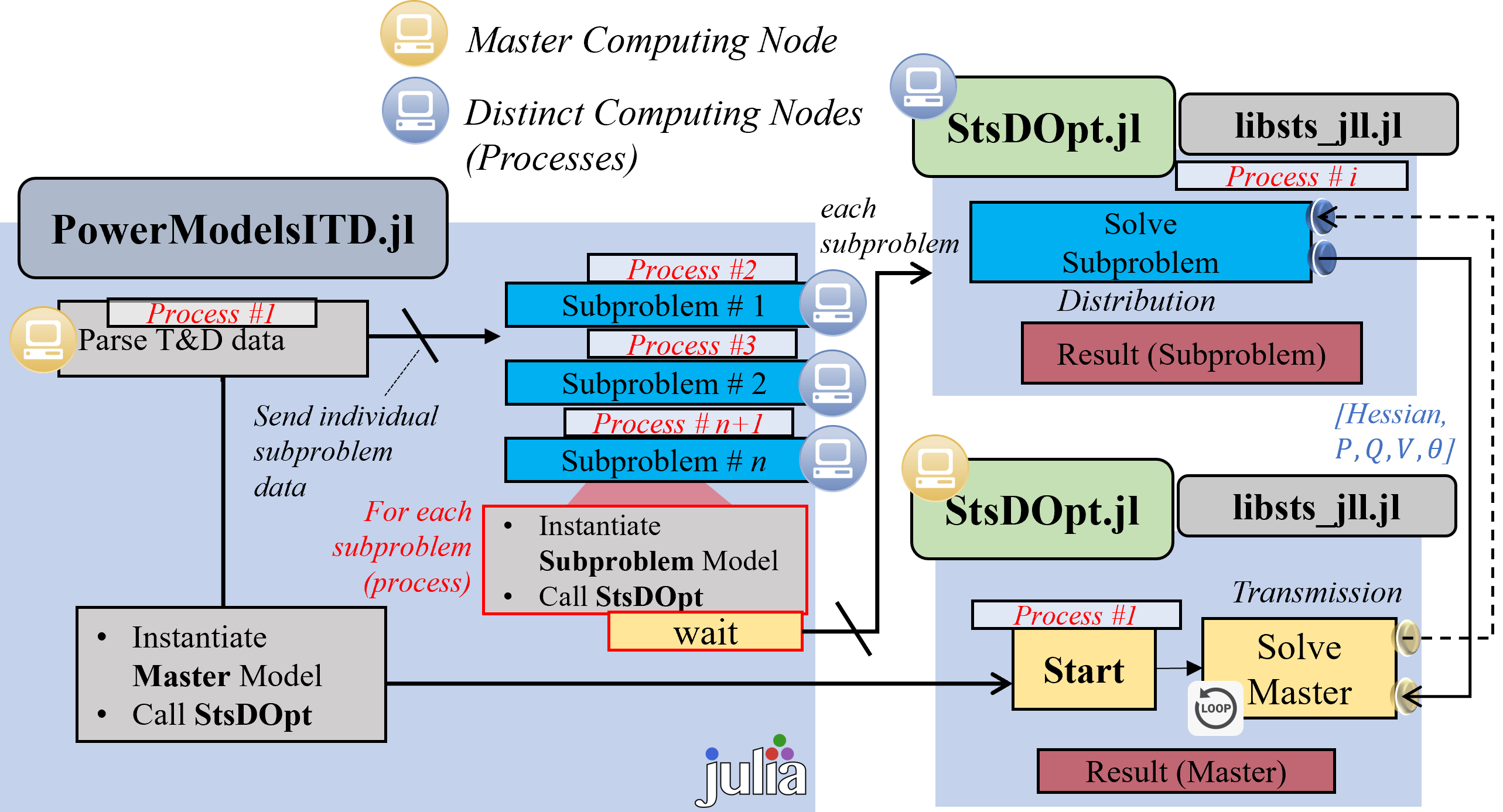}
\caption{\label{fig:multiprocess} Multiprocess version of StsDOpt.}
\end{figure}

This section discusses the implementation of the smoothed two-stage decomposition-based algorithm so that it can be seamlessly used within the Julia JuMP \cite{Lubin2023} environment, and more importantly within the PowerModelsITD optimization framework to enable large-scale T\&D AC-OPF optimization. Figs. \ref{fig:multithread} and \ref{fig:multiprocess} show diagrams that describe the 1) multithread and the 2) multiprocess versions of the StsDOpt solver developed based on smoothed two-stage decomposition algorithm presented in Section \ref{sect:sts}.

As observed in the diagrams, both the multithread and multiprocess versions of StsDOpt depend on two primary Julia packages, the `libsts\_jll.jl' and the `StsDOpt.jl'. These packages are essential for using the proposed solver across the Julia JuMP environment. JuMP is the optimization modeling language supported by the PowerModelsITD framework. A detailed description of these packages is presented in the next subsections, but before we dive into details, it is critical to discuss the differences between the multithread and multiprocess versions of the StsDOpt solver.

As seen in Fig. \ref{fig:multithread}, the multithread version of the StsDOpt solver is capable of performing the multithread operation when the `solving' loop is executing. The entire solving process, including the parsing, instantiation of the mathematical optimization models, and the execution of the solver, are performed in the same computing node (i.e., master node) and memory is shared between the the different threads allocated for solving the optimization problem. This version of the solver allows the user to allocate as many threads as possible (or available), depending on the number of cores, and each thread is capable of handling, in parallel, any subproblem assigned by the CPU scheduler. This is possible due to the fact that all the required memory is accessible by any single thread. The main advantage of the multithread version of the StsDOpt solver is that by using multiple threads, the time required to solve and instantiate large-scale T\&D AC-OPFs with many distribution systems can be significantly reduced compared to the monolithic integrated T\&D approach \cite{ospina2024modeling} because subproblems can be instantiated and solved in parallel.

The multiprocess version of the StsDOpt, as seen in Fig. \ref{fig:multiprocess}, executes the instantiation and solution processes associated with each subproblem on distinct computing nodes with distinct memory. The same instantiation and solution processes that exist in both the monolithic integrated T\&D AC-OPF and the multithread approaches need to be performed, with the difference that these can be executed in separate computing nodes or processes, and only subsets of the subproblems need to be loaded on the threads or processes. The master and subproblems are instantiated and initialized on separate processes, and the subproblem processes need to wait until the master problem gives them the signal to start solving their own optimization problems. A more complex coordination between the master and subproblems processes needs to be put in place, compared to the multithread version. Note that the current implementation of the multiprocess StsDOpt allows one compute node to handle multiple subproblems (as subprocesses). 

The multiprocess version of StsDOpt is primarily intended for cases when parallel resources are readily available or large-scale T\&D models are impossible to fit entirely in one single compute node due to lack of memory.  Using the multiprocess version of StsDOpt, users will be able to distribute the amount of memory and resources across the available computing nodes, significantly reducing the memory footprint of the problem. The multiprocess version of StsDOpt is a great fit for High-performance computing (HPC), where subproblems can be distributed throughout compute nodes in the cluster.

All the advantages and disadvantages between the different versions of the StsDOpt solver are discussed in Section \ref{sect:Results} using numerical results of large-scale T\&D AC-OPF problems. Now, let us discuss the packages developed that allow the StsDOpt solver to be used inside the PowerModelsITD framework.

\subsection{Implementation of the libsts\_jll.jl package}

The source code for the StsDOpt optimizer (i.e., the SQP master solver, the Ipopt-based subproblem solver, etc.), is written as a \textit{C++} software package. To enable its use for large-scale T\&D systems through the PowerModelsITD framework, we developed Julia JLL libraries that serve as wrappers to the callback functions required by the \textit{C++} software package. The Julia package that contains the pre-built libraries and executables is called libsts\_jll.jl, an auto-generated package constructed using BinaryBuilder.jl. Using BinaryBuilder, we are able to build shared libraries that can be used in different platforms such as Windows and GNU-Linux.

\subsection{Implementation of the StsDOpt.jl package}

The StsDOpt.jl package is the Julia wrapper for the StsDOpt solver that provides a high-level interface to the packages and libraries created by libsts\_jll.jl. This wrapper contains all the functions needed for users to ccall the StsDOpt solver and access its functionality. It enables the seamless interaction with the StsDOpt solver from within the Julia programming environment and contains the necessary orchestration processes (i.e., parallelization) required to handle the data exchange between the master problem and its corresponding subproblems. The StsDOpt.jl package is designed in a modular and flexible way that facilitates its use with any optimization tool that can decompose problems in a master-subproblem structure using the Julia JuMP modeling environment.

\subsection{Interfacing StsDOpt solver with PowerModelsITD}

The standardized development processes used for creating the StsDOpt.jl wrapper and the libsts\_jll.jl libraries makes the integration of the StsDOpt solver into the PowerModelsITD framework a seamless task. Only the parsing and instantiation procedures needed to be adapted since the master and subproblems are instantiated as separate mathematical optimization problems with their corresponding shared variables and boundary constraints, as opposed to the monolithic integrated T\&D AC-OPF approach used by default in PowerModelsITD \cite{ospina2024modeling}.

A common interface for the multithread and multiprocess versions of the solver is implemented, allowing for efficient exploitation of multithreaded, multi-core CPUs, and HPC clusters that enable faster solution times and lower memory footprints for large-scale T\&D AC-OPF instances. By utilizing multiple threads or processes, our implementation significantly improve the scalability of the StsDOpt solver. The solver is designed to be transparent and easy to use from the user's perspective, requiring minimal modifications to existing codebases that utilize PowerModelsITD.

\section{Numerical Results}\label{sect:Results}

\subsection{Description of Computational Setup}

The numerical experiments are conducted using a laptop computer and an HPC cluster. We employ a laptop with an Intel Xeon W-10855M CPU with 6 physical and 12 logical cores running at 2.80 GHz, and 80 GB of RAM. The operating system used is Ubuntu Linux 24.04, running under Windows Subsystem for Linux 2 (WSL2). The optimization solver Ipopt (version 3.14.16)\cite{wachter2006implementation} with the HSL MA27 linear solver (version 2023.11.17)\cite{hsl27}, is used for the solution of the integrated monolithic formulation used for benchmarking and the subproblems in the decomposition formulation. For the StsDOpt solver, the Ipopt source code required additional adaptations. The T\&D mathematical model is created using PowerModelsITD.jl (version 0.10.0) \cite{ospina2024modeling}, with an added direct interface to StsDOpt.jl (version 0.1.6), which, in turn, connects to libsts\_jll.jl (version 0.2.1).

The multiprocess experiments are run on a heterogeneous HPC cluster. Multiple compute nodes are utilized, with each node executing multiple subproblems managed using \textit{Slurm}. A maximum of 8 compute nodes are allocated with a maximum of 3 to 5 processes per node for the largest test cases. Each compute node has 20 physical cores running at 2.6 GHz, and 125 GB of RAM.

It is important to note that memory allocation for the multiprocess runs is primarily tracked on the master node, as data distributed across the cluster nodes is difficult to consolidate, and we are mostly interested in the memory footprint decrease in the master compute node compared to the multithread and the integrated monolithic approach.

\subsection{Description of Test Cases}

The performance of the integrated monolithic \cite{ospina2024modeling} approach and the multithread and multiprocess versions of the StsDOpt solver are compared using the T\&D test cases described next. All distribution system models are Kron reduced and all test cases are modeled using the nonlinear nonconvex AC polar formulation at both the transmission and distribution levels.

\begin{enumerate}
    \item \textbf{Case5-Case3 with 2 DGs}: A modified version of the PJM 5-bus system is used as the transmission system. The distribution system is connected as a new load added at bus \#5. The distribution system consists of a modified version of the IEEE 4 Node Test Feeder with increased loads and line capacities, and two distributed generators (DGs) rated at 1,000 kW and 200 kW added. 

    \item \textbf{Case5-Case3x2 with 2 DGs}: A modified version of the PJM 5-bus system is used as the transmission system, where distribution systems are connected at two new load buses, namely buses \#5 and \#6. The distribution systems consist of a modified version of the IEEE 4 Node Test Feeder with increased loads and line capacities. Two DGs rated at 1,000 kW and 200 kW are added to each of the feeders. One feeder is balanced, the other one is unbalanced.

    \item \textbf{Case118-Case3x5 with 2 DGs}: The IEEE 118-bus test system is used as the transmission system. Five distribution systems are connected at buses \#2, \#7, \#14, \#28, and \#44. The distribution systems consist of modified versions of the IEEE 4 Node Test Feeder with increased loads and line capacities, and one DG, with different power rating, added. The DGs' ratings are 600 kW, 400 kW, 600 kW, 200 kW, and 560 kW, respectively.

    \item \textbf{Case500-CaseLVx4}: The IEEE PGLib synthetic 500-bus test system is used as the transmission system. Four distribution systems are connected at buses \#323, \#337, \#114, and \#366. The distribution systems consist of modified versions of the IEEE LVTestCase (an European Low-Voltage test feeder) system with increased loads and line capacities, and three DGs (rated at 100 kW) added.
    
    \item \textbf{Case118-CaseR1R2GC}: The IEEE 118-bus test system is used as the transmission system. Ten distribution systems are connected at buses \#2, \#3, \#7, \#11, \#13, \#14, \#16, \#17, \#20, and \#21. The distribution systems consist of modified versions of the PNNL Taxonomy Feeders \cite{schneider2008modern} Regions 1, 2, and the General Circuit (GC). Multiple DERs, with different power ratings, are added to each of the feeders. The smallest feeder in this set has 96 nodes (the GC) and the biggest feeder has 3,491 nodes. The total number of nodes of the T\&D problem is 15,010 nodes.

    \item \textbf{Case118-CaseR3R4R5}: The IEEE 118-bus test system is used as the transmission system. 12 distribution systems are connected at buses \#23, \#28, \#29, \#33, \#35, \#39, \#41, \#43, \#50, \#51, \#52, and \#53. The distribution systems consist of modified versions of the PNNL Taxonomy Feeders \cite{schneider2008modern} Regions 3, 4, 5. Multiple DERs, with different power ratings, are added to each of the feeders. The smallest feeder in this set has 857 nodes and the biggest feeder has 11,367 nodes. The total number of nodes of the T\&D problem is 37,274 nodes.

    \item \textbf{Case118-CaseR1R2R3R4R5GC}: The IEEE 118-bus test system is used as the transmission system. 22 distribution systems are connected at buses \#2, \#3, \#7, \#11, \#13, \#14, \#16, \#17, \#20, \#21, \#23, \#28, \#29, \#33, \#35, \#39, \#41, \#43, \#50, \#51, \#52, and \#53. The distribution systems consist of modified versions of the PNNL Taxonomy Feeders \cite{schneider2008modern} Regions 1-5 and the GC. Multiple DERs, with different power ratings, are added to each of the feeders. The smallest feeder in this set has 96 nodes and the biggest feeder has 11,367 nodes. The total number of nodes of the T\&D problem is 52,166 nodes.

    \item \textbf{Case500-CaseR1R2GC}: The IEEE PGLib synthetic 500-bus test system is used as the transmission system. 10 distribution systems are connected at buses \#2, \#3, \#4, \#5, \#6, \#12, \#13, \#14, \#15, and \#16. The distribution systems consist of modified versions of the PNNL Taxonomy Feeders \cite{schneider2008modern} Regions 1, 2, and the GC. Multiple DERs, with different power ratings, are added to each of the feeders. The smallest feeder in this set has 96 nodes (the GC) and the biggest feeder has 3,491 nodes. The total number of nodes of the T\&D problem is 15,392 nodes.

    \item \textbf{Case500-CaseR3R4R5}: The IEEE PGLib synthetic 500-bus test system is used as the transmission system. 12 distribution systems are connected at buses \#17, \#19, \#21, \#23, \#24, \#29, \#55, \#56, \#58, \#60, \#62, and \#64. The distribution systems consist of modified versions of the PNNL Taxonomy Feeders \cite{schneider2008modern} Regions 3, 4, 5. Multiple DERs, with different power ratings, are added to each of the feeders. The smallest feeder in this set has 857 nodes and the biggest feeder has 11,367 nodes. The total number of nodes of the T\&D problem is 37,656 nodes.

    \item \textbf{Case500-CaseR1R2R3R4R5GC}: The IEEE 118-bus test system is used as the transmission system. 22 distribution systems are connected at buses \#2, \#3, \#4, \#5, \#6, \#12, \#13, \#14, \#15, \#16, \#17, \#19, \#21, \#23, \#24, \#29, \#55, \#56, \#58, \#60, \#62, and \#64. The distribution systems consist of modified versions of the PNNL Taxonomy Feeders \cite{schneider2008modern} Regions 1-5 and the GC. Multiple DERs, with different power ratings, are added to each of the feeders. The smallest feeder in this set has 96 nodes and the biggest feeder has 11,367 nodes. The total number of nodes of the T\&D problem is 52,548 nodes.
    
\end{enumerate}

The test cases were designed to gradually increase in size, allowing for comparisons that assess how the effectiveness of the StsDOpt solver improves as the problem size grows, and how memory allocation becomes critical for large-scale test cases.

\begin{table*}[htbp] \centering
\setlength{\tabcolsep}{1.8pt}
\caption{Comparisons of Integrated vs. multithread vs. multiprocess T\&D solvers on test cases.} 
\label{tab:numericalresults}
\begin{tabular}{|c|c|ccc|l|l|c|c|c|c|}
\hline
\textbf{\begin{tabular}[c]{@{}c@{}}Test \\ Cases\end{tabular}} & \textbf{\begin{tabular}[c]{@{}c@{}}Problem \\ Formulation\end{tabular}} & \multicolumn{3}{c|}{\textbf{\begin{tabular}[c]{@{}c@{}}\# of\\ Variables\end{tabular}}} & \multicolumn{1}{c|}{\textbf{\begin{tabular}[c]{@{}c@{}}Cost\\ (\$/hr)\end{tabular}}} & \multicolumn{1}{c|}{\textbf{\begin{tabular}[c]{@{}c@{}}Optimality \\ Gap\end{tabular}}} & \textbf{\begin{tabular}[c]{@{}c@{}}Wall-Clock \\ Solve Time (s)\end{tabular}} & \textbf{\begin{tabular}[c]{@{}c@{}}Total\\  Time (s)\end{tabular}} & \textbf{Iterations} & \textbf{\begin{tabular}[c]{@{}c@{}}Memory \\ Allocation\end{tabular}} \\ \hline \hline
\multirow{4}{*}{\begin{tabular}[c]{@{}c@{}}Case5-\\ Case3\end{tabular}} & \textbf{\begin{tabular}[c]{@{}c@{}}Integrated\end{tabular}} & \multicolumn{3}{c|}{214} & \$     17,947.59 &  & 0.03 & 0.04 & 27 & 11.9MiB \\ \cline{2-11} 
 & \textbf{} & \multicolumn{1}{c|}{\textbf{Master}} & \multicolumn{1}{c|}{\textbf{\begin{tabular}[c]{@{}c@{}}Sub\\ Min\end{tabular}}} & \textbf{\begin{tabular}[c]{@{}c@{}}Sub\\ Max\end{tabular}} &  &  &  &  &  &  \\ \cline{2-11} 
 & \textbf{Multithread} & \multicolumn{1}{c|}{54} & \multicolumn{1}{c|}{170} & 170 & \$     17,947.59 & 0.000\% & 0.26 & 0.28 & 30 & 12.3MiB \\ \cline{2-11} 
 & \textbf{\begin{tabular}[c]{@{}c@{}}Multiprocess\\  (x1 Processes)\end{tabular}} & \multicolumn{1}{c|}{54} & \multicolumn{1}{c|}{170} & 170 & \$     17,947.59 & 0.000\% & 0.27 & 0.56 & 30 & 12.5MiB \\ \hline \hline
\multirow{4}{*}{\begin{tabular}[c]{@{}c@{}}Case5-\\ Case3x2\end{tabular}} & \textbf{\begin{tabular}[c]{@{}c@{}}Integrated\end{tabular}} & \multicolumn{3}{c|}{384} & \$     18,092.98 &  & 0.07 & 0.09 & 26 & 21.7MiB \\ \cline{2-11} 
 & \textbf{} & \multicolumn{1}{c|}{\textbf{Master}} & \multicolumn{1}{c|}{\textbf{\begin{tabular}[c]{@{}c@{}}Sub\\ Min\end{tabular}}} & \textbf{\begin{tabular}[c]{@{}c@{}}Sub\\ Max\end{tabular}} &  &  &  &  &  &  \\ \cline{2-11} 
 & \textbf{Multithread} & \multicolumn{1}{c|}{64} & \multicolumn{1}{c|}{170} & 170 & \$     18,092.98 & 0.000\% & 0.27 & 0.29 & 30 & 22.3MiB \\ \cline{2-11} 
 & \textbf{\begin{tabular}[c]{@{}c@{}}Multiprocess\\  (x2 Processes)\end{tabular}} & \multicolumn{1}{c|}{64} & \multicolumn{1}{c|}{170} & 170 & \$     18,092.98 & 0.000\% & 0.37 & 0.38 & 30 & 6.05MiB \\ \hline \hline
\multirow{4}{*}{\begin{tabular}[c]{@{}c@{}}Case118-\\ Case3x5\end{tabular}} & \textbf{\begin{tabular}[c]{@{}c@{}}Integrated\end{tabular}} & \multicolumn{3}{c|}{1,858} & \$     94,567.56 &  & 0.21 & 0.27 & 26 & 94.1MiB \\ \cline{2-11} 
 & \textbf{} & \multicolumn{1}{c|}{\textbf{Master}} & \multicolumn{1}{c|}{\textbf{\begin{tabular}[c]{@{}c@{}}Sub\\ Min\end{tabular}}} & \textbf{\begin{tabular}[c]{@{}c@{}}Sub\\ Max\end{tabular}} &  &  &  &  &  &  \\ \cline{2-11} 
 & \textbf{Multithread} & \multicolumn{1}{c|}{1,108} & \multicolumn{1}{c|}{164} & 164 & \$     94,567.56 & 0.000\% & 3.78 & 3.86 & 49 & 99.5MiB \\ \cline{2-11} 
 & \textbf{\begin{tabular}[c]{@{}c@{}}Multiprocess\\  (x5 Processes)\end{tabular}} & \multicolumn{1}{c|}{1,108} & \multicolumn{1}{c|}{164} & 164 & \$     94,567.56 & 0.000\% & 4.35 & 4.40 & 49 & 59.6MiB \\ \hline \hline
\multirow{4}{*}{\begin{tabular}[c]{@{}c@{}}Case500-\\ CaseLVx4\end{tabular}} & \textbf{\begin{tabular}[c]{@{}c@{}}Integrated\end{tabular}} & \multicolumn{3}{c|}{70,094} & \$ 452,271.68 &  & 28.54 & 34.14 & 53 & 5.68GiB \\ \cline{2-11} 
 & \textbf{} & \multicolumn{1}{c|}{\textbf{Master}} & \multicolumn{1}{c|}{\textbf{\begin{tabular}[c]{@{}c@{}}Sub\\ Min\end{tabular}}} & \textbf{\begin{tabular}[c]{@{}c@{}}Sub\\ Max\end{tabular}} &  &  &  &  &  &  \\ \cline{2-11} 
 & \textbf{Multithread} & \multicolumn{1}{c|}{4,066} & \multicolumn{1}{c|}{16,520} & 16,520 & \$ 452,271.69 & 0.000\% & 54.64 & 63.47 & 53 & 5.89GiB \\ \cline{2-11} 
 & \textbf{\begin{tabular}[c]{@{}c@{}}Multiprocess\\  (x4 Processes)\end{tabular}} & \multicolumn{1}{c|}{4,066} & \multicolumn{1}{c|}{16,520} & 16,520 & \$ 452,271.69 & 0.000\% & 88.64 & 89.19 & 53 & 592MiB \\ \hline \hline
\multirow{4}{*}{\begin{tabular}[c]{@{}c@{}}Case118-\\ CaseR1R2GC\end{tabular}} & \textbf{\begin{tabular}[c]{@{}c@{}}Integrated\end{tabular}} & \multicolumn{3}{c|}{189,803} & \$     89,277.27 &  & 34.87 & 51.67 & 38 & 10.6GiB \\ \cline{2-11} 
 & \textbf{} & \multicolumn{1}{c|}{\textbf{Master}} & \multicolumn{1}{c|}{\textbf{\begin{tabular}[c]{@{}c@{}}Sub\\ Min\end{tabular}}} & \textbf{\begin{tabular}[c]{@{}c@{}}Sub\\ Max\end{tabular}} &  &  &  &  &  &  \\ \cline{2-11} 
 & \textbf{Multithread} & \multicolumn{1}{c|}{1,128} & \multicolumn{1}{c|}{1,118} & 45,634 & \$     89,277.40 & 0.000\% & 38.77 & 49.69 & 21 & 9.48GiB \\ \cline{2-11} 
 & \textbf{\begin{tabular}[c]{@{}c@{}}Multiprocess\\  (x10 Processes)\end{tabular}} & \multicolumn{1}{c|}{1,128} & \multicolumn{1}{c|}{1,118} & 45,634 & \$     89,277.40 & 0.000\% & 36.55 & 41.16 & 21 & 1.69GiB \\ \hline \hline
\multirow{4}{*}{\begin{tabular}[c]{@{}c@{}}Case118-\\ CaseR3R4R5\end{tabular}} & \textbf{\begin{tabular}[c]{@{}c@{}}Integrated\end{tabular}} & \multicolumn{3}{c|}{472,503} & \$     83,302.89 &  & 147.82 & 229.27 & 64 & 32.6GiB \\ \cline{2-11} 
 & \textbf{} & \multicolumn{1}{c|}{\textbf{Master}} & \multicolumn{1}{c|}{\textbf{\begin{tabular}[c]{@{}c@{}}Sub\\ Min\end{tabular}}} & \textbf{\begin{tabular}[c]{@{}c@{}}Sub\\ Max\end{tabular}} &  &  &  &  &  &  \\ \cline{2-11} 
 & \textbf{Multithread} & \multicolumn{1}{c|}{1,136} & \multicolumn{1}{c|}{9,770} & 140,198 & \$     83,305.85 & 0.004\% & 124.35 & 191.21 & 25 & 30.8GiB \\ \cline{2-11} 
 & \textbf{\begin{tabular}[c]{@{}c@{}}Multiprocess\\  (x12 Processes)\end{tabular}} & \multicolumn{1}{c|}{1,136} & \multicolumn{1}{c|}{9,770} & 140,198 & \$     83,305.85 & 0.004\% & 122.18 & 144.01 & 25 & 5.12GiB \\ \hline \hline
\multirow{4}{*}{\begin{tabular}[c]{@{}c@{}}Case118-\\ CaseR1R2R3R4R5GC\end{tabular}} & \textbf{\begin{tabular}[c]{@{}c@{}}Integrated\end{tabular}} & \multicolumn{3}{c|}{661,253} & \$     83,310.90 &  & 413.98 & 629.27 & 75 & 51.3GiB \\ \cline{2-11} 
 & \textbf{} & \multicolumn{1}{c|}{\textbf{Master}} & \multicolumn{1}{c|}{\textbf{\begin{tabular}[c]{@{}c@{}}Sub\\ Min\end{tabular}}} & \textbf{\begin{tabular}[c]{@{}c@{}}Sub\\ Max\end{tabular}} &  &  &  &  &  &  \\ \cline{2-11} 
 & \textbf{Multithread} & \multicolumn{1}{c|}{1,176} & \multicolumn{1}{c|}{1,118} & 140,198 & \$     83,315.21 & 0.005\% & 162.01 & 241.86 & 27 & 40.3GiB \\ \cline{2-11} 
 & \textbf{\begin{tabular}[c]{@{}c@{}}Multiprocess\\  (x22 Processes)\end{tabular}} & \multicolumn{1}{c|}{1,176} & \multicolumn{1}{c|}{1,118} & 140,198 & \$     83,315.21 & 0.005\% & 163.92 & 218.08 & 27 & 6.76GiB \\ \hline \hline
\multirow{4}{*}{\begin{tabular}[c]{@{}c@{}}Case500-\\ CaseR1R2GC\end{tabular}} & \textbf{\begin{tabular}[c]{@{}c@{}}Integrated\end{tabular}} & \multicolumn{3}{c|}{192,800} & \$ 410,093.09 &  & 40.25 & 57.81 & 40 & 10.7GiB \\ \cline{2-11} 
 & \textbf{} & \multicolumn{1}{c|}{\textbf{Master}} & \multicolumn{1}{c|}{\textbf{\begin{tabular}[c]{@{}c@{}}Sub\\ Min\end{tabular}}} & \textbf{\begin{tabular}[c]{@{}c@{}}Sub\\ Max\end{tabular}} &  &  &  &  &  &  \\ \cline{2-11} 
 & \textbf{Multithread} & \multicolumn{1}{c|}{4,090} & \multicolumn{1}{c|}{1,118} & 45,634 & \$ 410,094.24 & 0.000\% & 33.46 & 43.88 & 12 & 9.45GiB \\ \cline{2-11} 
 & \textbf{\begin{tabular}[c]{@{}c@{}}Multiprocess\\ (x10 Processes)\end{tabular}} & \multicolumn{1}{c|}{4,090} & \multicolumn{1}{c|}{1,118} & 45,634 & \$ 410,094.24 & 0.000\% & 28.82 & 33.65 & 12 & 1.82GiB \\ \hline \hline
\multirow{4}{*}{\begin{tabular}[c]{@{}c@{}}Case500-\\ CaseR3R4R5\end{tabular}} & \textbf{\begin{tabular}[c]{@{}c@{}}Integrated\end{tabular}} & \multicolumn{3}{c|}{475,500} & \$ 410,094.36 &  & 122.41 & 208.98 & 50 & 32.5GiB \\ \cline{2-11} 
 & \textbf{} & \multicolumn{1}{c|}{\textbf{Master}} & \multicolumn{1}{c|}{\textbf{\begin{tabular}[c]{@{}c@{}}Sub \\ Min\end{tabular}}} & \textbf{\begin{tabular}[c]{@{}c@{}}Sub\\ Max\end{tabular}} &  &  &  &  &  &  \\ \cline{2-11} 
 & \textbf{Multithread} & \multicolumn{1}{c|}{4,098} & \multicolumn{1}{c|}{9,770} & 140,198 & \$ 410,096.46 & 0.001\% & 100.32 & 147.10 & 18 & 25.5GiB \\ \cline{2-11} 
 & \textbf{\begin{tabular}[c]{@{}c@{}}Multiprocess\\  (x12 Processes)\end{tabular}} & \multicolumn{1}{c|}{4,098} & \multicolumn{1}{c|}{9,770} & 140,198 & \$ 410,096.46 & 0.001\% & 98.51 & 120.35 & 18 & 5.26GiB \\ \hline \hline
\multirow{4}{*}{\begin{tabular}[c]{@{}c@{}}Case500-\\ CaseR1R2R3R4R5GC\end{tabular}} & \textbf{\begin{tabular}[c]{@{}c@{}}Integrated\end{tabular}} & \multicolumn{3}{c|}{664,250} & \$ 410,102.68 &  & 225.19 & 422.28 & 52 & 50.8GiB \\ \cline{2-11} 
 & \textbf{} & \multicolumn{1}{c|}{\textbf{Master}} & \multicolumn{1}{c|}{\textbf{\begin{tabular}[c]{@{}c@{}}Sub \\ Min\end{tabular}}} & \textbf{\begin{tabular}[c]{@{}c@{}}Sub \\ Max\end{tabular}} &  &  &  &  &  &  \\ \cline{2-11} 
 & \textbf{Multithread} & \multicolumn{1}{c|}{4,138} & \multicolumn{1}{c|}{1,118} & 140,198 & \$ 410,105.58 & 0.001\% & 134.30 & 215.62 & 18 & 39.9GiB \\ \cline{2-11} 
 & \textbf{\begin{tabular}[c]{@{}c@{}}Multiprocess\\  (x22 Processes)\end{tabular}} & \multicolumn{1}{c|}{4,138} & \multicolumn{1}{c|}{1,118} & 140,198 & \$ 410,105.58 & 0.001\% & 117.27 & 146.44 & 18 & 6.91GiB \\ \hline
\end{tabular}
\end{table*}

\subsection{Numerical Results and Comparison Studies}
\label{sec:test1}
Using these test cases, we compare the performance and memory allocation footprint of the different T\&D approaches, i.e., the integrated monolithic T\&D approach against both the multithread and multiprocess StsDOpt solver (parallel approaches).

The metrics for comparison are the cost (in \$/hr), the optimality gap (in \%) - defined as:

\begin{equation}
    Opt_{gap} = \frac{\big(Obj. \ of \ Parall. \ - \ Obj. \ of \ Integr.\big)} {|Obj. \ of \ Integr.|} * 100
\end{equation}

, the wall-clock time to solve the problem (in seconds), the total time to instantiate and solve the problem (in seconds), the number of iterations of the respective solver, and the memory allocation (in Mebibyte (MiB) or Gibibyte (GiB)) as reported by Julia. ``Instantiating'' a problem refers to the process of building the JuMP models for the optimizers. 

Table \ref{tab:numericalresults} shows the numerical results of the test cases. As seen in the table, for the small test cases, such as \textit{Case5-Case3}, the integrated approach demonstrates superior performance with a solve time of 0.03 seconds and a total time of 0.04 seconds, significantly faster than both the multithread (0.26 seconds and 0.28 seconds, respectively) and  multiprocess (0.27 seconds and 0.56 seconds). This trend remains consistent as the problem size increases to \textit{Case5-Case3x2}, where the integrated approach maintains the fastest solve time of 0.07 seconds, whereas the multithread and multiprocess methods take 0.27 seconds and 0.37 seconds, respectively. However, the memory allocation for the multiprocess solver in \textit{Case5-Case3x2} is markedly lower (6.05 MiB) compared to both the integrated and multithread approaches (21.7 MiB and 22.3 MiB, respectively). This lower memory allocation requirement trend continues as the problem size increases.

For medium-size test cases, such as \textit{Case500-CaseLVx4}, the integrated approach's solve time of 28.54 seconds still outperforms the multithread (54.64 seconds) and the multiprocess (88.64 seconds) approaches. However, the multiprocess solver continues to exhibit a significantly reduced memory allocation of 592 MiB, when compared to the 5.68 GiB and 5.89 GiB used by the integrated and multithread approaches, respectively, due to its distributed nature. A similar pattern is observed in \textit{Case118-CaseR1R2GC}, where the integrated approach has a slightly faster solve time of 34.87 seconds when compared to the multithread approach's 38.77 seconds and the multiprocess approach's 36.55 seconds. Nonetheless, at this problem size and structure, the solve time starts equalizing and the total time is slightly lower for the decomposition approaches. Also, the memory allocation trend continues as the multiprocess solver only allocates 1.69 GiB of memory in the master computing node, significantly less when compared to the 9.48 GiB allocated by the multithread approach and the 10.6 GiB allocated by the integrated approach.

The performance disparity widens in the larger cases, such as \textit{Case118-CaseR3R4R5}, \textit{Case118-CaseR1R2R3R4R5GC}, \textit{Case500-CaseR1R2GC}, \textit{Case500-CaseR3R4R5}, and \textit{Case500-CaseR1R2R3R4R5GC}. For instance, in the largest case, \textit{Case500-CaseR1R2R3R4R5GC}, the integrated approach requires 225.19 seconds to solve the T\&D optimization problem, while the multithread and multiprocess approaches took 134.30 seconds and 117.27 seconds, respectively. A similar behavior is observed in the other large-scale test cases. The difference in memory allocation is also clearly seen here as the multiprocess solver’s memory usage was much lower at 6.91 GiB, in contrast to the 39.9 GiB allocated by the multithread approach and the 50.8 GiB required when using the integrated method.

In summary, the integrated approach consistently offered the fastest solve times across small- to medium-sized problems, where a small number of subproblems are being considered. However, as problem size, complexity, and the number of subproblems increases, the multiprocess method exhibits a significant advantage in memory efficiency and solution time, making it a viable option for scenarios where memory constraints are critical and multiple computing nodes are available. The multithread approach generally performs well, exhibiting substantial advantages in solution and instantiation times compared to the integrated approach, suggesting that its application is advantageous in specific cases when only one computing node is available and the problem consists of multiple distribution systems (subproblems) that can be solved concurrently using the presented decomposition approach.

\vspace{-4mm}
\subsection{Scalability Studies}

The ultimate goal of our research is the optimization of very large-scale power systems consisting of hundreds of distribution systems connected to the same transmission system network. However, currently, there are no large-scale test cases, realistic or synthetic, that can be used for evaluation. Instead, to get an idea of the scalability of the proposed decomposition approach, we consider a scalable network for which the mathematical structure, in terms of size and the number of subproblems, is similar to realistic large-scale instances.

The network has the IEEE 118-bus test system as the transmission system, to which multiple copies of the PNNL Taxonomy R2-12.47-3 distribution system (comprising 3,491 nodes) are added at one bus. The size of the network can be adjusted by changing the number of copies.  We ran the proposed methods for instances with 1, 2, 10, 20, and 40 copies, resulting in optimization problems with 44, 179, 87, 305, 432, 313, 863, 573, and 1,726,093 variables, respectively. The three approaches are compared in terms of solve time (in seconds) and total time (in seconds).

\begin{figure}[h]
\centering
\includegraphics[width = 0.45\textwidth]{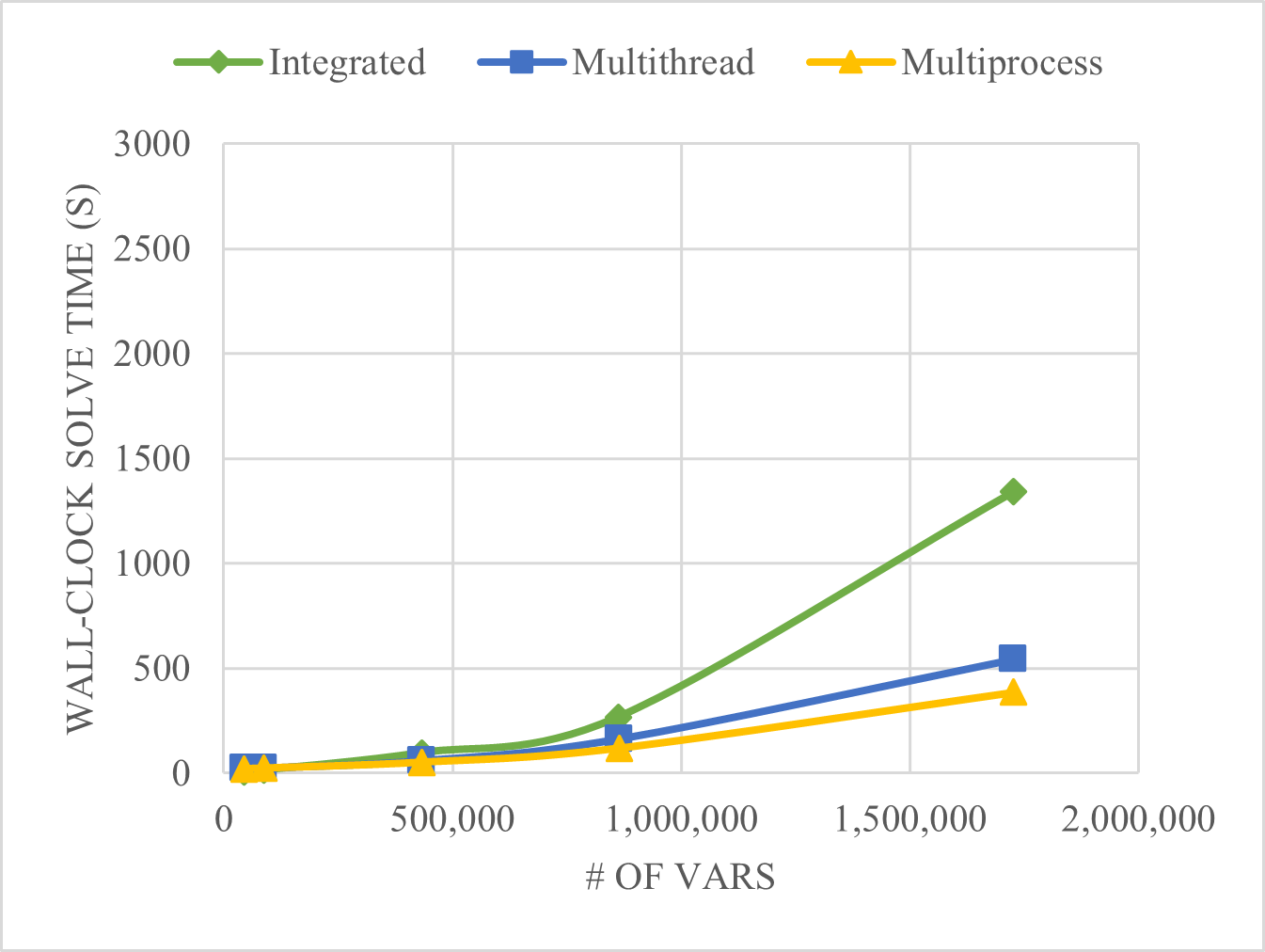}
\caption{\label{fig:varswallclock} Number of variables in T\&D AC-OPF problem vs. Wall-clock solve time in seconds for scalability case.}
\end{figure}

\begin{figure}[h]
\centering
\includegraphics[width = 0.45\textwidth]{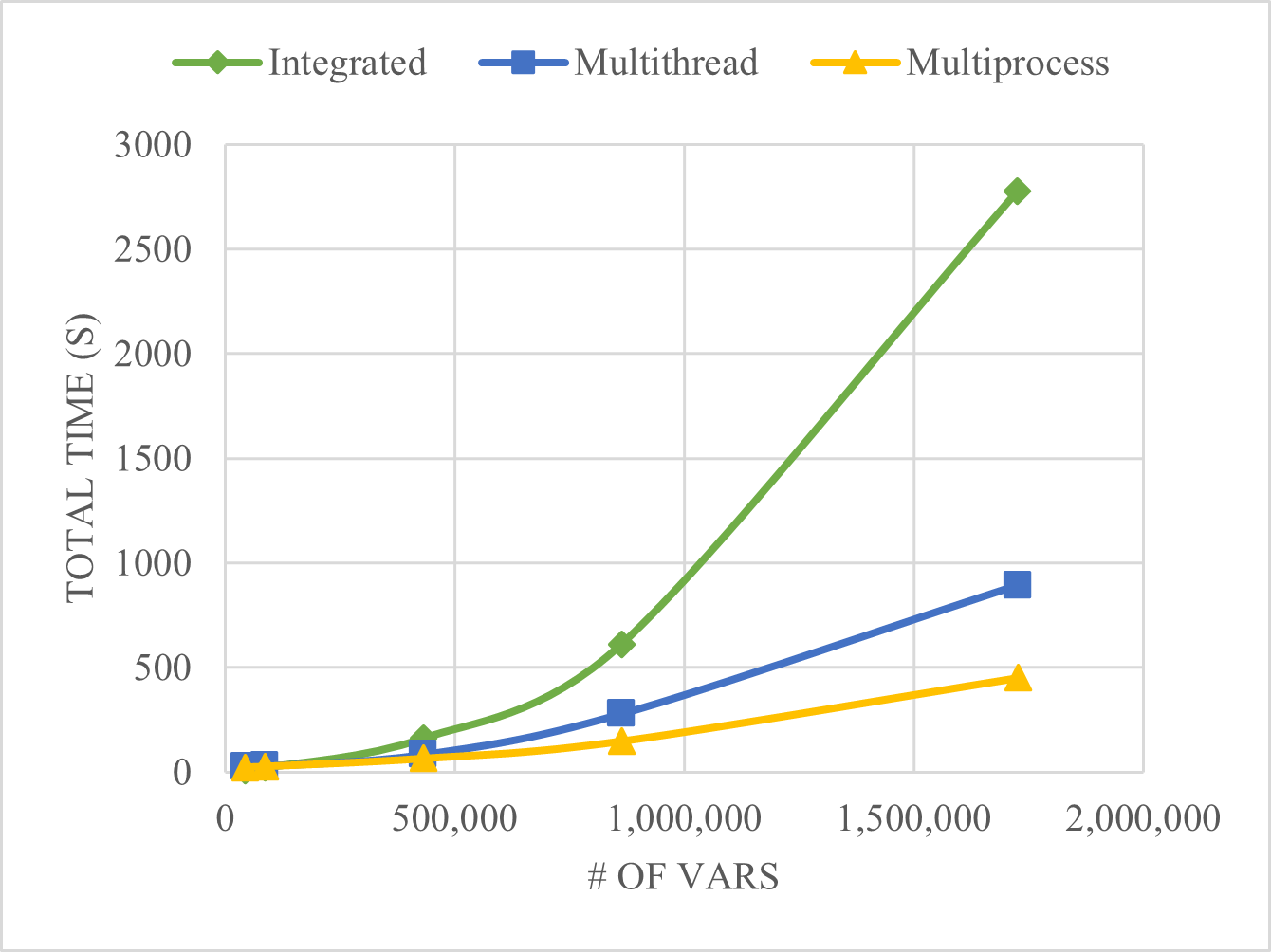}
\caption{\label{fig:varsinstantclock} Number of variables in T\&D AC-OPF problem vs. Total time in seconds for scalability case.}
\end{figure}

Figs.~\ref{fig:varswallclock} and \ref{fig:varsinstantclock} plot the solution time and total time against the number of optimization variables. As the number of variables and subproblems increases, the advantage of using the proposed decomposition approach, whether in multithread or multiprocess mode, becomes increasingly evident. These trends are consistent with the observations discussed in Section~\ref{sec:test1}.

\vspace{-5mm}
\section{Conclusions}
\label{sect:conclusions}

This paper presents a novel smoothed two-stage decomposition optimizer (StsDOpt) designed to efficiently solve large-scale transmission and distribution (T\&D) AC optimal power flow (AC-OPF) problems. The proposed method offers several advantages over existing solutions, including its ability to handle nonlinear and nonconvex models without relying on linearization or relaxation, as well as its ability to handle accurate modeling of unbalanced and nonlinear distribution systems. 

Numerical results demonstrate the effectiveness of the proposed methodology in solving large-scale, nonlinear, nonconvex T\&D AC-OPF problems. Computational experiments were conducted across a variety of test cases, scaling up to 22 distribution systems and approximately 650k nodes. These experiments demonstrated that the proposed optimizer remains scalable and efficient, even under memory constraints, as evidenced by its performance in multithread and multiprocess environments. Compared to the integrated T\&D approach, StsDOpt consistently achieved reduced computational times maintaining solution accuracy in large-scale test cases.

The scalability of the StsDOpt solver was further validated through studies that measured its performance as the number of distribution systems increased. These studies revealed that StsDOpt outperforms traditional methods, maintaining solution quality while reducing computation time, making it particularly suitable for large-scale AC-OPF problems. This `scalability’ is crucial as power grids become increasingly complex with higher levels of DER penetration, necessitating optimization techniques capable of managing the evolving complexity and scale of T\&D systems.

The proposed solver is integrated into the PowerModelsITD optimization framework, allowing for thorough testing and verification of the results, ensuring the accuracy and reliability of the presented results.

Overall, the proposed method presents a promising solution for solving large-scale AC-OPF problems with high-levels of DER penetration, offering significant advancements in computational efficiency and scalability. Future work will focus on improving solver convergence, incorporating additional AC-OPF formulations such as AC rectangular and IV rectangular, and developing larger test cases to fully explore the limits of the proposed methodology.
\vspace{-2mm}
\section*{Acknowledgments}
\label{sect:acknowledgements}
This work was performed with the support of the U.S. Department of Energy (DOE) Office of Electricity (OE) Advanced Grid Modeling (AGM) Research Program under program manager Ali Ghassemian. This research used resources provided by the Darwin testbed at Los Alamos National Laboratory (LANL), funded by the Computational Systems and Software Environments subprogram of LANL's Advanced Simulation and Computing program. The research work conducted at LANL is done under the auspices of the National Nuclear Security Administration of the U.S. Department of Energy under Contract No. 89233218CNA000001.
\vspace{-2mm}
\bibliography{bibfile}
\bibliographystyle{IEEEtran}

\end{document}